\documentclass[12pt]{amsart}
\usepackage[english]{babel}
\usepackage{amsmath}
\usepackage{amsthm}
\usepackage{amssymb}
\usepackage{mathrsfs}
\usepackage{enumerate}
\usepackage{graphicx}
\usepackage{epstopdf}
\usepackage{float}
\usepackage{caption}
\usepackage[notcite, final, notref]{showkeys}
\usepackage{dsfont}
\usepackage{array}
\usepackage[dvips]{color}
\newtheorem{theorem}{Theorem}[section]

\newtheorem{definition}[theorem]{Definition}

\theoremstyle{definition}
\newtheorem{remark}[theorem]{Remark}
\newtheorem{example}[theorem]{Example}

\def\E{\bf E}

\usepackage{xcolor}

\title[]{Sparse Approximation to the Dirac-$\delta$ Distribution}

\author[W. Qu]{Wei Qu}
\address{(WQ)School of Mathematical Sciences\\
 Beijing Normal University \\
 China}
\email{quwei2019@bnu.edu.cn}

\author[T. Qian]{Tao Qian*}
\address{(TQ)Macau Center for Mathematical Sciences\\
Macau University of Science and Technology\\
Macau}
\email{tqian@must.edu.mo}
\thanks{*Corresponding author\\
Funded by The Science and Technology Development Fund, Macau SAR (File no. 0123/2018/A3)}

\author[G.-T. Deng]{Guan-Tie Deng}
\address{(GTD)School of Mathematical Sciences\\
 Beijing Normal University\\
 China}
\email{96022@bnu.edu.cn}

\def\R{\bf R}

\begin{document}
\maketitle
\begin{abstract}
The Dirac-$\delta$ distribution may be realized through sequences of convlutions, the latter being also regarded as approximation to the identity. The present study proposes the so called pre-orthogonal adaptive Fourier decomposition (POAFD) method to realize fast approximation to the identity. The type of sparse representation method has potential applications in signal and image analysis, as well as in system identification.
\end{abstract}

\bigskip


\bigskip


\date{today}
\tableofcontents
\def\ux{\underline{x}}
\def\uy{\underline{y}}
\def\uxi{\underline{\xi}}

\def\HH{$\mathcal{H}$-$H_K$}
\def\bD{\bf D}

\section{Introduction}
\setcounter{equation}{0}
The most common examples of approximation to the identity include those of the convolution integral type by using the Poisson kernel, the heat kernel, and some more general convolution kernels satisfying certain normalization conditions (\cite{SW}). In the series form we have Poisson summation etc. From these classical examples one can observe that a signal may be well approximated by a finite linear combination of the convolution kernel of the context. In this study we develop an approximation theory of such type.
\def\H{\mathcal{H}}
The approximation can be associated with an axiomatic or text-book formulation, that we call $\H$-$H_K$ formulation (\cite{Q2020,Saitoh1,Saitoh2}), of Hilbert space with a linear operator defined through an inner product kernel.
We give a quick revision on this formulation. Let
$\mathcal{H}$ be a general Hilbert space with inner product
$\langle\cdot ,\cdot\rangle_{\mathcal H}.$
Let ${\bf E},$ the set of parameters, be a set of numbers, or a set of vectors, whose components are real or complex numbers. $\E$ is assumed to be an open set with respect to the usual topology of ${\bf R}^n$ or ${\bf C}^n.$
Let every $p\in \E$ be
\def\C{\bf C}
associated with an element $h_p\in \mathcal{H}$ that gives rise to
 a linear
operator $L: \mathcal{H}\to {\bf C}^{\bf E},$ the latter being the set of all functions from $\E$ to $\C.$
 \begin{eqnarray}\label{operatorL}
 L(f)(p)=\langle f,h_p\rangle_{\mathcal H}.\end{eqnarray}
 We also write $F(p)=Lf(p)$ and denote by $R(L)$ the function set of the images $F\in {\bf C}^{\bf E}$ of $f\in \mathcal{H}$ under the mapping $L.$
Let $N(L)$ be the null space of $L$ defined as
\[ N(L)=\{f\in \mathcal{H}\ |\ L(f)=0\}.\] It is easy to show that
$N(L)$ is a closed set in $\mathcal{H}.$
 There thus exists the orthogonal complement of $N(L)$ in $\mathcal{H}$ denoted $N(L)^\perp,$ and
\[ \mathcal{H}=N(L)\oplus N(L)^\perp.\]
Accordingly, each $f\in \mathcal{H}$ can be uniquely written as
\[ f=f^-+f^+,\]
where $f^-\in N(L), f^+\in N(L)^\perp.$ In the set-mapping notation, there holds
$L(N(L)^\perp)=R(L).$
Denote the orthogonal projection operator from
$\mathcal{H}$ to $N(L)^\perp$ by $P: P(f)=f^+.$ We introduce a new Hilbert space
structure, $H_K,$  on the function set $R(L),$ as follows:
The induced norm of $F=L(f)$ in the range set $R(L)$ is defined as
$$\|F\|_{H_K}\triangleq\|Pf\|_{\mathcal{H}}.$$
The norm definition induces an inner product in $R(L)$ denoted
$\langle \cdot,\cdot\rangle_{H_K}.$   The new Hilbert space $H_K,$ coinciding with
$R(L)$ in the set-theoretic sense, is
isometric with $N(L)^\perp$ through the mapping $L.$
 In such notation the function
$K(p,q)$ defined
\[ K(q,p)=\langle h_q,h_p\rangle_{\mathcal{H}}\] is, in fact, the reproducing kernel of $H_K,$ and hence, the latter is a reproducing kernel Hilbert space. For a proof of this, see \cite{Q2020} or \cite{Saitoh1}.
While the $\H$-$H_K$ formulation makes it convenient to study linear operator theory in Hilbert spaces in general this paper will concentrate in the particular case where the class of functions $\{h_p\}_{p\in \E}$ is a dense subspace of $\H.$ In the case the null space $N(L)$ is trivial, containing only the zero function. In fact,
\[ \langle f,h_p\rangle=0 \qquad \forall p\in \E\]
if and only if
\def\D{\bf D}
\[ L(f)(p)=0\qquad \forall p\in \E,\]
and thus $N(L)=\{0\}$ and $N(L)^\perp=\H={\rm span}\{h_p\}_{p\in \E}.$ In the case it would be very beneficial and instructive that although $\H$ is not a RKHS but $H_K$ is, the latter being isometric with the former under the mapping $L.$
The density of $\{h_p\}_{p\in \E}$ in $\H$ amounts that $\H$ is a Hilbert space with a dictionary $\{h_p/\|h_p\|\}_{p\in \E}.$ In such sense any separable Hilbert space, therefore has a dictionary, is equivalent with a RKHS. The latter enjoys useful properties that offer more technical methods in dealing with separable Hilbert spaces. The best example of $\H$-$H_K$ structure is $\H=L^2(\partial \D),$ the $L^2$ space on the unit circle, and ${\E}={\D}, h_p(e^{it})=\frac{1}{1-\overline{p}e^{it}}, p\in {\D}.$ In the case $N(L)^\perp=H_+^2(\partial {\D})=R(L), N(L)=H_-^2(\partial {\D}),$ being respectively the boundary Hardy spaces inside and outside the unit circle. It is as if the $\H$-$H_K$ formulation is specially made for this and the other classical Hardy spaces situation, but actually not, for the structure is possessed by all linear operator in Hilbert spaces induced by a kernel with a parameter, and in particular includes all linear differential and integral operators. Paper \cite{Q2020} initiates the sparse solutions methodology to basic problems of operators defined through inner product kernels within the $\H$-$H_K$ formulation.

The goal of this article is to introduce a particular sparse representation method in general Hilbert spaces possessing a dictionary with the so called \emph{boundary vanishing condition} (BVC). The sparse representation is called pre-orthogonal adaptive Fourier decomposition, or POAFD in brief (see \S 2). POAFD is a generalization of the so called adaptive Fourier decomposition, or AFD, originally developed for the classical complex Hardy $H^2$ spaces. The AFD for one complex variable right fits into the delicate frame work of the Beurling-Lax Theorem involving Blaschke products (\cite{QWa}). Some engineering applications of AFD may be found, for instance, in \cite{Mi1,Mi2,MQW,LZQ,WQLG}. Some generalizations of AFD to higher dimensions are successful \cite{QSW,ACQS1,ACQS2,Q2D}. Because of lack of Blaschke product or Takenaka-Malmquist system in context, generalizations of AFD to domains other than the classical types, or to multi-dimensions or analytic function spaces other than the Hardy type, however, are difficult or impossible.  POAFD, with general applicability, reduces to AFD in the classical Hardy space case, being of the ultimate optimality among various types of greedy algorithms (see \cite{Te} and \cite{MZ}): The POAFD maximal selection is the greediest among all the one-step-optimal selections. POAFD is, in particular, supported by repeating selection of the parameters, involving,  when necessary, Gram-Schmidt orthogonalization of directional derivatives of the dictionary elements.

 If in an $\mathcal{H}$-$H_K$ formulation the function set $\{h_p/\|h_p\|_{\mathcal{H}}\}_{p\in {\bf E}}$ is a dictionary of the underlying space $\mathcal{H},$ then there exist two equivalent approaches to construct the POAFD type sparse representation in $\mathcal{H}.$
 One is a directly application of POAFD in $\mathcal{H}$ just by using the dictionary properties. The other is to perform the sparse representation in $H_K,$ which has the advantage as a RKHS, in which we have the convenience to normalize the kernels and to prove BVC. After getting a sparse series expansion in $H_K$ we convert back the obtained expansion to $\mathcal{H}$ through the
  isometric mapping $L^{-1}.$
The purpose of this study is to develop a general sparse representation methodology for the Dirac-$\delta$ distribution with the understanding and help from the point of view of $\H$-$H_K$ formulation. In particular, the RKHS approach brings in delicate analysis and helps in getting better understanding to the subject.

In \S2 we give a detailed description of the POAFD method. In \S 3 we develop the convolution type sparse representation of the identity in the underlying space ${\R}^d$ using POAFD, including the Poisson and heat kernels and the general convolution kernel satisfying non-degenerate and the usual decaying rate at the infinity. In \S 4 we develop, as a bounded $\E$ case,  the spherical Poisson sparse approximation to the identity. Having given detailed description of the POAFD method in general Hilbert spaces with a dictionary satisfying BVC in \S 2, what we do in \S 3 and \S 4, as the main body of the paper, are verifications of its applicability to the most common and yet important models, i.e., the Poisson and the heat kernels, convolution kernels in general, as well as the spherical Poisson kernel case. The verifications are proceeded under the frame work of the $\H$-$H_K$ formulation. \S 5 contains two illustrative examples on, respectively, the spherical Poisson POAFD on the sphere and heat kernel POAFD in ${\R}^d.$

\section{POAFD in Hilbert Space with a Dictionary}
The basic idea and the related concepts, including POAFD maximal selection principle, boundary vanishing condition and multiple kernels, first appeared in \cite{Q2D}. The formulation of the method, including terminology in use, has been revised and improved, and unified through a sequence of related studies. At the beginning POAFD is designed for sparse representation of images defined on rectangles, being topologically identical with the $2$-torus. Later this method is extended to spaces of analytic functions other than Hardy spaces (\cite{MQ1,MQ2,QD1,QD2,CQT}). We now introduce some related concepts in Hilbert spaces with a dictionary.

In the $\mathcal{H}$-$H_K$ formulation $H_K$ is a RKHS with the kernel function $K(q,p)=K_q(p)=
\langle h_q,h_p\rangle_{\mathcal{H}}.$ Since $\langle F,K_q\rangle=0$ for all $q\in {\bf E}$ implies $F=0,$ we know that the function set $\{K_q\}_{q\in {\bf E}}$ is dense in $H_K.$
\begin{definition} A subset $\mathcal{E}$ of a general Hilbert space ${H}$ is said to be a \emph{dictionary} if $\|E\|=1$ for $E\in \mathcal{E},$ and
$\overline{\rm span}\{E
 :\ E\in \mathcal{E}\}={H}.$
\end{definition}
With the notation of the last section the normalized reproducing kernels $E_q=K_q/\|K_q\|, q\in \E,$ constitute a dictionary of $H_K.$ On the $\H$ space side in any case the functions $h_p/\|h_p\|, p\in {\E},$ constitute a dictionary of $N(L)^\perp;$ and, if $\{h_p\}_{p\in \E}$ is dense in $\mathcal{H},$ then the functions $h_p/\|h_p\|, q\in {\E},$ constitute a dictionary of $\mathcal{H}.$

The POAFD method is available for all Hilbert spaces that has a dictionary, regardless whether the dictionary is from a reproducing kernel or not. In below we sometimes borrow the notation $K_q, q\in \E,$ not assuming their reproducing property but only assuming density of $\{K_q\}_{q\in \E}.$  The normalized form $E_q=K_q/\|K_q\|$ is used only when involving the so called boundary vanishing condition (BVC, see below).

Before we introduce the maximal selection
principle of POAFD we need to introduce two concepts: \emph{boundary vanishing condition} (BVC) and \emph{multiple reproducing kernel}. With BVC we need to make some convention when $\E$ is an unbounded set in its underlying space, say ${\R}^{d+1}.$ This will be the case when we discuss the Poisson and the heat kernels in the following sections, in which $\E$ is the upper-half space of ${\R}^{d+1}.$ In the case we add one more point, $\infty,$ to the whole space ${\R}^{d+1}.$ We make $\infty$ to be a new boundary point of $\E$ by modifying the topology of ${\R}^{d+1}$ through introducing an open neighborhood system of $\infty:$ A set $O$ is said to be an open neighborhood of $\infty$ if and only if the complement of $O,$
$O^c$ is a compact set of ${\R}^{d+1}.$ That is, we use the compactification of ${\R}^{d+1}$ with respect to the added point $\infty.$  We denote by $\partial^*\E$ the set $\partial \E\cup \{\infty\},$ which is the set of boundary points of $\E$ in the new topological space ${\R}^{d+1}\cup\{\infty\},$ where $\partial \E$ is the set of all finite boundary points of $\E.$  As a consequence, an open neighbourhood of $\partial^* \E$ is the union of an open neighbourhood of the set $\partial \E$ and an open neighborhood of $\infty.$ Since under the one-point-compactification topology the space ${\R}^{d+1}\cup\{\infty\}$ is compact, its closed subset $\E\cup\partial^*\E$ is also compact.

If $\E$ is a bounded open set in ${\R}^{d+1},$ such as when we discuss spherical Poisson kernel approximation, then we do not have to do anything with the original topology. In such case we are with the convention $\partial^*\E=\partial \E.$ Boundary vanishing condition (BVC) in both the bounded and unbounded $\E$ cases are stated as
\begin{definition} Let ${H}$ be a Hilbert space with a dictionary $\mathcal{E}=\{E_q\}_{q\in \E}.$ If for any $f\in \mathcal{H}$ and any $q_k\to \partial^* \E$, in the one-point-compactification topology if necessary, there holds
 \[ \lim_{k\to \infty} |\langle f,E_{q_k}\rangle|=0,\]
then we say that ${H}$ together with $\mathcal{E}$ satisfy BVC.\end{definition}

If ${H}$ is a Hilbert space with a dictionary  $E_q, q\in \E,$ satisfying BVC, then a compact argument will lead to the conclusion for any $f\in H$ that there exists a selection of $\tilde{q}\in \E$ such that
 \[ \tilde{q}=\arg \sup\{ |\langle f,E_{{q}}\rangle|^2\ :\ q\in \E\}.\]
 We note that the Hilbert space $H$ in the definition can be, with regards to the  $\H$-$H_K$ model, two cases: $H=\H$ or $H=H_K.$ With the case $H=\H$ we refer to the dictionary $\{h_p/\|h_p\|\}_{p\in \E},$ while in the second case $H=H_K$ we refer to the
 collection of the normalized reproducing kernels $E_q=K_q/\|K_q\|, q\in \E.$

  Many RKHSs, including the classical Hardy spaces, Bergman and weighted Bergman spaces, satisfy BVC. On the other hand, there exist RKHSs whose normalized reproducing kernels constituting a dictionary that does not satisfy BVC ([QD1]).

Next we define multiple kernels. Let $(q_1,\cdots,q_n)$ be an $n$-tuple of parameters in ${\bf E}.$ The set ${\bf E}$ may be a region
 in the complex plane, or one in ${\bf R}^d,$ or even in ${\bf C}^d.$ In the ${\bf R}^d$ case let
\begin{eqnarray}\label{multiple} \tilde{K}_{q_n}(p)=\left[\left(\frac{\partial}{\partial q_{\vec{\theta}}}\right)^{j(n)-1}K_q\right]_{q=q_n}(p),\end{eqnarray}
where $j(n)$ is the number of repeating times of the parameter $q_n$ in the  $n$-tuple
$(q_1,\cdots,q_n),$ in the case $\frac{\partial}{\partial q_{\vec{\theta}}}=\vec{\theta}\cdot\nabla,$ being the directional derivative in the direction $\vec{\theta}=(\cos\theta_1,\cdots,\cos\theta_d).$  If ${\E}\subset {\bf C}^d,$ the concept is similarly defined. For the case ${\E}\subset {\bf C}$ the directional derivative is simply replaced by $e^{i\theta}\frac{\partial}{\partial z}.$ With such notation, if there is no repeating, that is $q_k\ne q_n$ for all $k<n,$ then $j(n)=1,$ and $\tilde{K}_{q_n}={K}_{q_n}.$   The kernel $\tilde{K}_{q_n}(p)$ is called the \emph{multiple kernel} associated with the $n$-tuple $(q_1,\cdots,q_n).$  Associated with an infinite sequence $(q_1,\cdots,q_n,\cdots ),$ there is, in such way, an infinite sequence of multiple
kernels $\tilde{K}_{q_n}, n=1,\cdots.$  In this paper we tacitly assume that all the involved directional derivatives of the dictionary elements of any order belong to the underlying Hilbert space $\mathcal{H}.$ The derivatives of one or higher orders occur during the optimization process through maximal selections of the parameter (\cite{Q20}, or its close English version \cite{CQT}, also \cite{QSW,Q2D,QD1}). We also call
a finite or infinite sequence of multiple kernels as a consecutive multiple kernel sequence, for if
$\tilde{K}_{q_n}$ involves the $j(n)-1$ order derivative, then all the proceeding $k$-derivative kernels for $k<j(n)-1,$ at the same point $q_n,$  should also have appeared before $\tilde{K}_{q_n}$ in the sequence.
 Together with BVC, the multiple kernels enable realization of the maximal
 selection principle in the following pre-orthogonal optimal process. Suppose we already have an $(n-1)$-tuple
 $(q_1,\cdots,q_{n-1}),$ allowing repetitions, and accordingly have
  an associated $(n-1)$-tuple of consecutive multiple kernels $\{\tilde{K}_{q_1},\cdots,\tilde{K}_{q_{n-1}}\}.$
 By performing the Gram-Schmidt (G-S) orthonormalization process we have
 an $(n-1)$-orthonormal tuple
 $\{B_1,\cdots,B_{n-1}\}$ that is equivalent, in the linear span sense, with $\{\tilde{K}_{q_1},\cdots,\tilde{K}_{q_{n-1}}\}.$
 The decisive role of BVC and multiple kernels is as follows: For any $g\in \mathcal{H},$
 one can find a $q_n\in \E$ such that
 \begin{eqnarray}\label{sup} q_n={\rm arg}\sup\{|\langle g,B_n^q\rangle|\ :\ q\in {\bf E}\},\end{eqnarray}
where $B_n^q$ is such that $\{B_1,\cdots,B_{n-1},B_n^q\}$ is the G-S orthonormalization
of $\{\tilde{K}_{q_1},\cdots,\tilde{K}_{q_{n-1}}, \tilde{K}_{q}\},$ and $B_n^q$ is precisely given by
\begin{eqnarray}\label{GS}
 B_n^q=\frac{\tilde{K}_q-\sum_{k=1}^{n-1}\langle \tilde{K}_q,B_k\rangle B_k}
 {\sqrt{\|\tilde{K}_q\|^2-\sum_{k=1}^{n-1}|\langle \tilde{K}_q,B_k\rangle|^2}}.\end{eqnarray}
 In the POAFD algorithm we will use this for $g=g_n,$ being the $n$-th standard remainder
\[ g_n=f-\sum_{k=1}^{n-1}\langle f
,B_k\rangle B_k.\]
 In such way we consecutively extract out the maximal energy portion from the standard remainders. At the step-by-step optimal selection category POAFD is, indeed,
  the greediest optimization strategy, that is guaranteed by BVC and the concept multiple kernels.
  The evolution of the idea and the exposition of POAFD can be found in the literature \cite{QWa,QSW,Q2D,Q20,QD1,CQT}.

 \begin{remark} If ${H}$ does not have a dictionary satisfying BVC then even with multiple kernels one cannot perform POAFD. However, from the definition of supreme, for any $\rho\in (0,1)$ and any mutually distinguished $q_1,\cdots,q_{n-1},$ there exists $q_n$ different from the preceding $q_k, k=1,\cdots,n-1,$ such that
\begin{eqnarray}\label{weak q} |\langle g,B_n^{q_n}\rangle|\ge \rho
\sup\{\langle g,B_n^q\rangle| \ :\ q\in {\bf \E}, q\ne q_1,\cdots,q_{n-1}\}.\end{eqnarray}
The algorithm for consecutively finding $q_n, n=1,\cdots,n,\cdots,$  to satisfy
(\ref{weak q}) is called \emph{Weak Pre-orthogonal
adaptive Fourier decomposition} (Weak-POAFD). Practically we often adopt the Weak-POAFD maximal principle, as, in the weak manner, we can at every step select a parameter different from what have been chosen in the previous steps. Theoretically, however, we are more interested in
the case where existence of the exact maximizers $q_n$ to (\ref{sup}) can be guaranteed. In the classical Hardy space case POAFD is equivalent with AFD using TM systems. Indeed, it can be proved that TM systems are not only orthonormal by themselves, but also are G-S orthonormalizations of the multiple Szeg\"o kernels of the context.\end{remark}

\def\bR{\bf R}

By using
POAFD one can prove that the $n$-th standard remainder of a POAFD is dominated by the magnitude $M/\sqrt{n}$ if the expanded function $f$ belongs to the space
\[ {H}^M=\{ f\ |\ f\in {H}, \exists \ q_k, d_k\ {\rm such \ that}\ f=\sum_{k=1}^\infty d_kE_{q_k} \ {\rm with} \ \sum_{k=1}^\infty |d_k|\leq M\}\]
 (see \cite{Q2D,QD1}).

 We remark that the above convergence rate estimation is promising as there is no smoothness condition imposed to the expanded function. With concrete examples usually much more rapid convergence are observed. As having in mind, the POAFD method is to be promoted with
the $\mathcal{H}$-$H_K$ formulation in numerical solutions of integral and differential equations (see \cite{Q2020}). In the present paper we only explore its impact with spars representation of the Dirac-$\delta$ distribution (\cite{SW}).

\def\F{\bf F}

\def\bR{\bf R}
\def\bC{\bf C}
\section{Sparse Approximation of the Convolution Type}

\subsection{Sparse Poisson Kernel Approximation}

 It is well known that Poisson integrals approximate the boundary data function. In this section we will develop sparse approximation by linear combinations of parameterized Poisson kernels.

The Poisson kernel context fits well with the  $\mathcal{H}$-$H_K$ formulation. We let $\mathcal{H}=L^2({\bf R}^d).$  Set
$${\bf E}=\{ p\in {\bf R}^{d+1}_+\ |\ p=t+\underline{x}, t>0, \underline{x}=(x_1,\cdots x_d)\}.$$ For $p=t+\underline{x},$
 let
\[ h_p(\uy)=P_{t+\underline{x}}(\underline{y})\triangleq c_d\frac{t}{|p-\uy|^{d+1}}=c_d
\frac{t}{(t^2+|\underline{x}-\underline{y}|^2)^{\frac{d+1}{2}}},\quad d\ge 1,\]
where $c_d=\frac{\Gamma [(d+1)/2]}{\pi^{(d+1)/2}}.$ We note that $h_p(\uy)$ is the evaluation  at the point $\ux-\uy$ of the $L^1$-$t$-dilation of the function $\phi(\ux)=c_d\frac{1}{(1+|\ux|^2)^{(d+1)/2}},$ where $c_d$ is the normalization constant under which
the integral of $\phi$ over ${\bR}^d$ is identical with $1.$ In below when we discuss the Poisson kernel on the unit sphere and the heat kernel in ${\bR}^d$ we use $c_d$ for the same normalization purpose, whose values then vary from context to context.

 The operator $L$ and its images $Lf, f\in L^2({\bR}^d),$ are given by
 \[ u(t+\underline{x})=Lf(t+\underline{x})=\langle f,h_{t+\underline{x}}\rangle_{L^2({\bR}^d)}.\]
  In the $\mathcal{H}$-$H_K$ formulation the range $R(L)$ consists of the Poisson
 integrals of the boundary data $f\in L^2({\bf R}^d).$ Now we show that $\{h_p\}_{p\in\E}$ is dense. It suffices to show that if $f\in L^2({\bR}^d)$ and
 $\langle f,h_p\rangle =0$ for all $p,$ then $f=0.$ It is a result of harmonic analysis that, in both the $L^2({\bR}^d)$-norm and pointwise sense,
 \[ \lim_{t\to 0+}u(t+\underline{x})=\lim_{t\to 0+}\langle f,h_p\rangle =f(\ux)=0,\quad {\rm a.e.}\]
 In the $\H$-$H_K$ formulation we have $N(L)=\{0\}$ and $N(L)^\perp=L^2({\bR}^d).$
 On the other hand $R(L)=H_K$ is, under the mapping $L$, isometric with $N(L)^\perp=L^2({\bR}^d).$ In particular, $L(h_q)=K_q.$ Density of $K_q$ in $H_K$ implies density of $h_p$ in $\H=L^2({\bR}^d).$

 Harmonic analysis knowledge has given a characterization of the space $H_K.$ In fact, $H_K$ coincides, together with its norm, with the harmonic Hardy space on the upper-half space ${\bR}^{d+1}_+:$
 \[ h^2({\bR}^{d+1}_+)=\{ u: {\bR}^{d+1}_+ \to {\bR}\ :\ \triangle_{{\bR}^{d+1}_+} u=0, \|u\|_{h^2({\bR}^{d+1}_+)}^2=\sup_{t>0}\int_{{\bR}^d}
 |u(t+\underline{x})|^2d\underline{x}<\infty\}.\]
By denoting
 $f(\underline{x})=u(0+\underline{x}),$ we have
 \begin{eqnarray*} \|u\|_{H_K}^2\overset{{\mathcal{H}{\mbox -}H_K}}{=}\|f\|^2_{L^2({\bR}^d)}
 \overset{{NBL}}{=}\|u(0+\cdot)\|^2_{L^2({\bR}^d)}\overset{{h^2{\mbox -}{\rm Theory}}}{=}
 \sup_{t>0}\int_{{\bR}^d}
 |u(t+\underline{x})|^2d\underline{x}.\end{eqnarray*}

 The reproducing kernel of the space $H_K$ is computed as, for
 $p=t+\underline{x}, q=t_1+\underline{x}_1,$
\begin{eqnarray}\label{thanks1}
K(q,p)&=&\langle P_{t_1+\underline{x}_1},P_{t+\underline{x}}\rangle_{L^2({\bR}^d)}\nonumber\\
&=&\int_{{\bR}^d} P_{t_1+\underline{x}_1}
(\underline{\xi})P_{t+\underline{x}}(\underline{\xi})d\underline{\xi}\nonumber \\
&=&c_d^2
\int_{{\bR}^d} \frac{t_1}{(t_1^2+|\underline{x}_1-\underline{\xi}|^2)^{\frac{d+1}{2}}}
\frac{t}{(t^2+|\underline{x}-\underline{\xi}|^2)^{\frac{d+1}{2}}}
d\underline{\xi}\nonumber \\
&=& P_{(t_1+t)+\underline{x}_1}(\underline{x})
\end{eqnarray}
 This last equality relation is due to the uniqueness of the solution
 of the Dirichelet problem:
 \[ \triangle u=0;  \quad {\rm and }\quad u(0+\underline{x})=P_{t_1+\underline{x}_1}(\underline{x}).\]
 Indeed, on one hand, the first three expressions of the above equality chain all mean that the left-hand-side is the Poisson integral of the boundary data $P_{t_1+\underline{x}_1}(\underline{\cdot}),$ and thus is harmonic in $t+\underline{x}.$ On the other hand, the function $P_{(t_1+t)+\underline{x}_1}(\underline{x})$, being harmonic in the variable $t+\underline{x}$, has the boundary limit function $P_{t_1+\underline{x}_1}(\underline{\cdot}).$  Therefore, these two harmonic functions have to be the same.

The above deduction also concludes the relation
\[\langle P_{t_1+\underline{x}_1},P_{t+\underline{x}}\rangle_{L^2({\bR}^d)}=
P_{(t_1+t)+(\underline{x}_1-\underline{x})}(0),\]
regarded as the semigroup property of the Poisson kernel.
  The reproducing property
 of $K_q$ is an immediate consequence of the $\H$-$H_K$ formulation: For $u\in H_K, q=t_1+\underline{x}_1,$
 \begin{eqnarray*}
 \langle u,K_q\rangle_{H_K}&=&\langle u(0+\underline{\cdot}),P_{t_1+\underline{x}_1}(\underline{\cdot})\rangle_{L^2({\bR}^d)}
 =u(t_1+\underline{x}_1).\end{eqnarray*}
 For a general
   $K_{q}, q=t+\underline{x}, t>0,$ its norm is computed, from the semi-group property (\ref{thanks1}),
   \[  \|K_q\|^2_{H_K}=\langle K_q,K_q\rangle_{H_K}=K(q,q)=P_{2t}(0)
   =\frac{c_d}{(2t)^d}.\]
   The norm-one normalization of $K_q$ is thus
  \[ E_q=\frac{K_q}{\|K_q\|}=\left(\frac{(2t)^{d}}{c_d}\right)^{1/2}K_q.\]
   Next we verify that BVC holds in this Poisson context, i.e.,
 \begin{eqnarray}\label{BVC}
  \lim_{q\to \partial \E}|\langle u,E_q\rangle_{H_K}|=0,\end{eqnarray}
 where $u$ is any function in $H_K=h^2({\bR}^{d+1}_+).$  We first have, by using the reproducing kernel property, for $q=t+\underline{x},$
 \begin{eqnarray}\label{from}
 \langle u,E_q\rangle_{H_K}=c_d't^{d/2}u(t+\underline{x}).
 \end{eqnarray}
Due to density of the parameterized Poisson kernels
in $H_K,$ the verification of BVC is reduced to verifying (\ref{BVC}) for each parameterized
 reproducing kernel $u(p)=K_{t_1+\underline{x}_1}.$ From (\ref{from}) we have
 \begin{eqnarray}\label{quantity}\langle K_{t_1+\underline{x}_1},E_q\rangle_{H_K}
 =c_dt^{d/2}P_{(t_1+t)+\underline{x}_1}(\underline{x})=
 c_d't^{d/2}\frac{t+t_1}{[(t+t_1)^2+
 |\underline{x}-\underline{x}_1|^2]^{(d+1)/2}}.
 \end{eqnarray}

 The limiting process $q\to \partial^* {\bR}^{d+1}_+,$ based on the one-point-compactification topology, amounts to, alternatively,
 $t\to 0$ or $t^2+\underline{x}^2\to
 \infty.$ For any fixed $\underline{x}_1$ and $t_1>0,$ regardless the positions of $\underline{x},$  we have
 \[
 t^{d/2}\frac{t+t_1}{[(t+t_1)^2+
 |\underline{x}-\underline{x}_1|^2]^{(d+1)/2}}\leq
 t^{d/2}\frac{1}{(t+t_1)^d}\to 0,\]
 as $t\to 0$ ($d\ge 1$). So,  uniformly in $\underline{x},$ as $t\to 0,$ the quantity in (\ref{quantity}) tends to zero.

 Let, for the fixed $t_1$ and $\underline{x}_1,$  $R=\sqrt{t^2+|\underline{x}|^2}>4|\underline{x}_1|+2t_1+1.$ We divide the argument into the two cases: (1) $0<t<R/2;$ and
 (2)\ $t\ge R/2.$ In case (1), $|\underline{x}|>R/2$  and  hence $(|\underline{x}-\underline{x}_1|)^2\ge (|\underline{x}|/4)^2.$  Hence, by ignoring the constant,
 \[ (t)^{d/2}\frac{t+t_1}{[(t+t_1)^2+
 |\underline{x}-\underline{x}_1|^2]^{(d+1)/2}}\leq
 (R/2)^{d/2}\frac{R/2+R/2}{[(t+t_1)^2+|\underline{x}/4|)^2]^{(d+1)/2}}\leq \frac{c_d}{R^{d/2}}\to 0,\]
 as $R\to \infty$ ($d\ge 1$).

 In case (2), $t\ge R/2$ implies
 \[ t^{d/2}\frac{t+t_1}
 {[(t+t_1)^2+|\underline{x}-\underline{x}_1|^2]^{(d+1)/2}}\leq
 t^{d/2}\frac{1}{(t+t_1)^{n}}\leq \frac{c_d'''}{R^{d/2}}\to 0,\]
 as $R\to \infty$ ($d\ge 1$). Thus,  uniformly in $t,$ when $\underline{x}$ tends to infinity the quantity in (\ref{quantity}) tends to zero.
\def\bS{{\bf S}^{d-1}}
\def\bB{{\bf B}^d}

BVC is thus proved. The maximizer of (\ref{sup}) is attainable in $\E.$ Hence POAFD can be proceeded.

 \subsection{Sparse Heat (Gaussian) Kernel Approximation}

 It is known that the integral transformation
 \begin{eqnarray}\label{write} Lf(t+\underline{x})=\frac{1}{(4\pi t)^{d/2}}\int_{{\bf R}^d}f(\underline{y})
 e^{-\frac{|\underline{x}-\underline{y}|}{4t}}d\underline{y}, \quad n\ge 1,\end{eqnarray}
 gives rise to the unique solution $u(t+\underline{x})$ of the initial value problem for the heat equation
 \[ \frac{\partial u}{\partial t}=\triangle u, \quad u(0+\underline{x})=f(\underline{x}), \quad f\in L^2({\bf R}^d),\]
 where $\triangle$ is the Laplacian for $x_1,\cdots,x_d.$
 This well fits into the $\mathcal{H}$-$H_K$ formulation with
 $\mathcal{H}=L^2({\bf R}^d), q=t+\underline{x}\in {\E}={\bf R}^{d+1}_+,$ and
 \begin{eqnarray}\label{definition} h_q(\underline{y})=\frac{1}{(4\pi t)^{d/2}}e^{-\frac{|\underline{x}-\underline{y}|}{4t}}.\end{eqnarray}
 The classical heat kernel analysis asserts that $\{h_p\}_{p\in {\bf R}^{d+1}_+}$ is dense in $L^2({\bR}^d).$
 The space $N(L)$ is hence the trivial subspace consisting of only the zero function,
 while the space $N(L)^\perp$ coincides with $L^2({\bR}^d).$ In the $\mathcal{H}$-$H_K$ formulation
 the space $H_K$ is as the range set $R(L)$ equipped with the induced norm from their $L^2$-boundary data. Associated with the heat kernel, the space $H_K$ may be characterized, like the defining conditions for $h^2({\bR}^{d+1}_+)$ in the Poisson kernel case, by using a quantization condition plus a condition such as a solution of an linear differential equation. Here we do not explore the details for they are not used for the main purpose of the study. The space
 $L^2({\bR}^d)$ coincides with the
 non-tangential boundary limits (NBL) of the Gauss-Weierstrass integral $u=Lf:$
 \[ \lim_{t\to 0}u(t+\underline{x})=f(\underline{x}),\quad f\in L^2({\bR}^d),\]
 in both the $L^2({\bR}^d)$- and the a.e. pointwise sense. For this reason we denote
 $f(\underline{x})=u(0+\underline{x}),$ and, as in the Poisson integral case, have the relations
 \begin{eqnarray*} \|u\|_{H_K}^2\overset{{\mathcal{H}{\mbox -}H_K}}{=}\|f\|^2_{L^2({\bR}^d)}
 \overset{{NBL}}{=}\|u(0+\cdot)\|^2_{L^2({\bR}^d)}.\end{eqnarray*}
 For the heat kernel case the space $H_K$ has a similar characterization as for the Poisson kernel case using the harmonic $h^2$ space. For our approximation purpose we only deduce the reproducing kernel. With
 $q=t+\underline{x}$ and $p=s+\underline{y}:$
 \begin{eqnarray}\label{infact}
 K(q,p)&=&\langle h_q,h_p\rangle_{L^2({\bf R}^d)}\nonumber \\
 &=&\frac{1}{(4\pi)^d(ts)^{d/2}}\int_{{\bf R}^d}e^{\frac{-|\underline{x}-\underline{\xi}|^2}{4t}}
 e^{\frac{-|\underline{y}-\underline{\xi}|^2}{4s}}d\xi.
 \end{eqnarray}
 We claim that the last integral representing $K(q,p)$ is equal to
 \begin{eqnarray}\label{new}\frac{1}{(4\pi)^d(ts)^{d/2}}\int_{{\bf R}^d}e^{\frac{-|\underline{x}-\underline{\xi}|^2}{4t}}
 e^{\frac{-|\underline{y}-\underline{\xi}|^2}{4s}}d\xi=
 \frac{1}{(4\pi(t+s))^{d/2}}e^{\frac{-|\underline{x}-\underline{y}|^2}{4(t+s)}},
 \end{eqnarray}
 and therefore, according to (\ref{definition}),
 \begin{eqnarray}\label{proved} K(q,p)=h_{(t+s)+\underline{x}}(\underline{y})=h_{(t+s)+\underline{y}}(\underline{x})
 =h_{(t+s)+(\underline{x}-\underline{y})}(0).
 \end{eqnarray}
 The proof of the identity (\ref{new}) uses the same idea as that in the Poisson kernel case, but involves the heat initial value problem
 \[  \frac{\partial u}{\partial t}=\triangle u, \quad \lim_{s\to 0+}u(s+\underline{y})=h_{(t+0)+\underline{x}}(\underline{y}).\]
The relation
\[\langle h_{t+\underline{x}},h_{s+\underline{y}}\rangle_{L^2({\bR}^d)}=
P_{(t+s)+(\underline{x}-\underline{y})}(0),\]
is regarded as the semigroup property of the heat kernel.

To proceed with heat kernel sparse representation using POAFD we first verify BVC.
The norm of the kernel $K_q$ is computed through
\[ \|K_q\|_{H_K}^2=\langle K_q,K_q\rangle_{H_K}=h_{2t+\underline{0}}(\underline{0})=
\frac{1}{(8\pi t)^{d/2}}.\]
Hence, the normalized reproducing kernel $K_q$
for $q=t+\underline{x}, p=s+\underline{y}$ becomes
\begin{eqnarray}\label{becomes} E_q(p)=(8\pi t)^{d/4}h_{(t+s)+\underline{x}}(\underline{y})=
\frac{(8\pi t)^{d/4}}{(4\pi(t+s))^{d/2}}
e^{-\frac{|\underline{x}-\underline{y}|^2}{4(t+s)}}.\end{eqnarray}
We are to show BVC, i.e.,
 \begin{eqnarray}\label{H1}
  \lim_{{\bf R}^{d+1}_+\ni q\to \partial^* {\bf R}^{d+1}_+}|\langle u,E_q\rangle_{H_K}|=0,
  \end{eqnarray}
 where $u$ is any function in $H_K.$
 Due to the density of the heat kernel in $L^2({\bR}^d)$ (\cite{SW}) the verification of BVC is reduced to only for an arbitrary but fixed reproducing kernel.
 As in the Poisson kernel case we use the one-point-compactification topology.  We need to show, for any but fixed $p=s+\underline{y},$ under the process $q=t+\underline{x}\to \partial^* {\bf R}^{d+1}_+$ we have
 \[ \lim_{q\to \partial^* {\bf R}^{d+1}_+} |\langle K_p,E_q\rangle_{H_K}|=0.\]
From the previous computation we have
 \begin{eqnarray}\label{quantity2}\langle K_p,E_q\rangle_{H_K}
 =E_q(p)=(8\pi t)^{d/4}h_{(t+s)+\underline{x}}(\underline{y}).
 \end{eqnarray}
 Write, for a constant $C$, the right hand side of (\ref{becomes}) as
 \begin{eqnarray}\label{accor} C\left(\frac{t}{t+s}\right)^{d/4}
 \left(\frac{1}{t+s}\right)^{d/4}e^{-\frac{|\underline{x}-\underline{y}|^2}{4(t+s)}}.\end{eqnarray}
 From (\ref{accor}) for the fixed $p=s+\underline{y}$ there exists a constant $C,$ depending only on the dimension $d,$ such that for all  $\underline{x}$ uniformly,
 \[ (8\pi t)^{d/4}h_{(t+s)+\underline{x}}(\underline{y})\leq C\left(\frac{t}{s+t}\right)^{d/4}\frac{1}{ s^{d/4}}\to 0, \quad {\rm as}\ t\to 0.\]
 Next we analyze the process $R=\sqrt{t^2+|\underline{x}|^2}\to \infty.$ Let
 $R\ge R_0=4|\underline{y}|+2s+1.$ We divide the argument into two cases: (1) $0<t<R/2;$ and
 (2)\ $t\ge R/2.$ In case (1), $|\underline{x}|>R/2,$  and  hence $(|\underline{x}-\underline{y}|)^2\ge (R/4)^2.$ In the case we have $t+s<R.$
 Thus, from (\ref{accor}),
 \[ (8\pi t)^{d/4}h_{(t+s)+\underline{x}}(\underline{y})\leq C\frac{1}{ s^{d/4}}e^{-\frac{(R/4)^2}{R}}\to 0, \quad {\rm as}\ R\to \infty.\]
 In case (2) we have $t\ge R/2.$ Through a brutal estimation based on (\ref{accor}) we have
 \[ (8\pi t)^{d/4}h_{(t+s)+\underline{x}}(\underline{y})\leq C\left(\frac{1}{t+s}\right)^{d/4}\leq C\frac{1}{R^{d/4}}\to 0, \quad {\rm as} \ R\to \infty.\]
 Thus BVC is proved and POAFD can be performed.

 \subsection{Sparse Approximation for the General Convolution Case}

Let $\phi\in L^1({\bR}^d)\cap L^2({\bR}^d)$ with
 \[ \int_{{\bR}^d}\phi (\ux)d\ux =1.\]
 We further assume the following conventional condition:
\begin{eqnarray}\label{decay} \psi({\ux})=\sup_{|{\uy}|\ge |{\ux}|}|\phi ({\uy})|\leq \frac{C}{(1+|{\ux}|^2)^{\frac{d+\delta}{2}}}, \quad \delta>0.\end{eqnarray}
Under these conditions we have the approximation to identity property in the $L^2$-sense (\cite{SW}):
For $\phi_t({\ux})=\frac{1}{t^d}\phi (\frac{\ux}{t}),$
\begin{eqnarray}\label{due to}\lim_{t\to 0}f\ast\phi_t=f,\ {\rm in\ } \ L^2,\
{\rm and\ pointwise\ as\ well}\quad \lim_{t\to 0}f\ast\phi_t({\ux})=f({\ux}), \quad {\rm a.e.}\end{eqnarray}
Let $\mathcal{H}={\bR}^d, {\E}={\bR}^{d+1}_+.$ Write, as before, $q=t+\ux, p=s+\uy, t,s>0,$ and $\ux, \uy\in {\bR}^d.$ In this general context $h_p(\uy)$ has the form
\[ h_p(\uy)=\frac{1}{t^d}\phi (\frac{\ux-\uy}{t}).\]
In the $\mathcal{H}$-$H_K$ formulation the $H_K$ space is
\[ H_K=\{u:{\bR}^{d+1}_+\to {\bR} \ :\ u(p)=\langle f,h_p\rangle_{L^2({\bR}^d)}\}.\]
In harmonic analysis characterizations of the space $H_K$ may involve a quantization condition such as $L^2$-boundedness of the non-tangential maximal function together with some none quantization but characterising property of the convolution kernel $\phi$ itself.  The details, however, are not needed in sparse representation study of this paper.
The reproducing kernel is computed
\begin{eqnarray}\label{kernel form}
 K_q(p)=K(q,p)&=&\langle h_q,h_p\rangle_{L^2({\bR}^d)}\nonumber\\
 &=&\frac{1}{(ts)^d}\int_{L^2({\bR}^d)}\phi (\frac{\ux-\uxi}{t})\phi (\frac{\uy-\uxi}{s})d\xi.\end{eqnarray}
 For $u=Lf$ the reproducing kernel property is automatic:
 \[ \langle u,h_p\rangle_{H_K}=\langle f(\cdot ),h_{p}(\cdot )\rangle_{L^2({\bR}^d)}=u(p).\]

The Poisson and heat kernels are particular cases of the above convolution form formulation, except that the heat kernel uses the replacement $\sqrt{t}$ for $t.$  In the Poisson and the heat kernel cases the dilations $\phi_t(\ux)$ and $\phi_{\sqrt{t}}(\ux),$ respectively, satisfy certain partial differential equations, and the convolutions against the boundary data give rise to the unique solution of the corresponding boundary value problem. In such cases one can prove certain semi-group property and have the formula $K_q(p)=h_{t+s+\ux}(\uy).$ In the general cases the kernel $K_q$ given by the integral (\ref{kernel form}) does not have semi-group property, not have a closed form either. In such case by using the integral formula (\ref{kernel form}) and the decaying property (\ref{decay}) one is able to prove BVC for $\delta\ge 1$.

 In fact, as before, a density argument reduces the verification for a parameterized reproducing kernel. In such case one can enlarge the last integral in (\ref{kernel form}) to get a Poisson kernel domination, and then get the limit as for the Poisson kernel case. Precisely, for a fixed $p=s+\uy$  and $ q=t+\ux\to\partial {\bR}^{d+1}_+,$ we may show
\[ |\langle K_p,E_q \rangle_{H_K}|\leq Ct^{d/2}P_{(s+t)+\underline{y}}(\underline{x})\to 0.\]
For the cases $0<\delta<1$ we do not know the answer.

We finally note that POAFD can be applied to either of the two contexts: the $\mathcal{H}$ context or the $H_K$ context. On the $H_K$ context we have reproducing kernel properties to use, that is very convenient especially when the kernel has an explicit formula.

\section{Poisson Kernel Sparse Approximation on Spheres}
 Next we set $\mathcal{H}=L^2(\bS),$ the Hilbert space of the square integrable functions on the $(d-1)$-dimensional unit sphere centered at the origin, and ${\E}=\bB,$ the $d$-dimensional open unit ball centered at the origin, where $d\ge 2.$
For a point $q\in \bB,$ the function $h_q$ in the context is the Poisson kernel of the ball: with $q=rt, r=|q|, t,s\in \bS,$
\begin{eqnarray}\label{q} h_q(s)=P_q(s)=c_d\frac{1-r^2}{|q-s|^d}.\end{eqnarray}
 The operator $L$ and its images $Lf, f\in L^2(\bS),$ are given by
 \[ u(q)=Lf(q)=\langle f,h_{q}\rangle_{L^2(\bS)},\]
 where the inner product of $L^2(\bS)$ is
 \[ \langle f,g\rangle_{L^2(\bS)}=\int_{\bS}f(s)g(s)d\sigma(s),\]
 where $d\sigma(s)$ is the normalized Lebesgue measure on the sphere. 
 The range $R(L)$ is the harmonic Hardy space on the unit ball $\bB:$
 \[ h^2({\bB})=\{ u: \bB \to {\bR}\ :\ \triangle u=0,\ \sup_{0\leq r<1}\int_{\bS}
 |u(rs)|^2d\sigma(s)<\infty\}.\]
 Due to the density of $\{h_q\}_{q\in \bB}$ in $L^2(\bS)$
 the space $N(L)$ is trivial consisting of only the zero function,
 while the space $N(L)^\perp$ is identical with $L^2(\bS).$
 The space $H_K=h^2(\bB),$ being the range set $R(L)$ equipped with the norm induced from their
 non-tangential boundary limits: For $u=Lf,$
 \[  \|u\|_{H_K}\triangleq \|f\|_{L^2(\bS)}.\]
 It is well known knowledge that
 \[ \lim_{r\to 1}u(rt)=f(t)\]
 in both the $L^2$- and in the a.e. pointwise sense. For this reason we write
 $f(t)=u(t)$ under the correspondence $u=Lf.$ We also have the relations
 \begin{eqnarray*} \|u\|_{H_K}^2\overset{{\mathcal{H}{\mbox -}H_K}}{=}\|f\|^2_{L^2(\bS)}
 \overset{{NBL}}{=}\|u(\cdot)\|^2_{L^2(\bS)}\overset{{h^2{\mbox -}{\rm Theory}}}{=}
 \sup_{0\leq r<1}\int_{\bS}
 |u(rt)|^2d\sigma(t).\end{eqnarray*}

 The reproducing kernel of the space $H_K$ is computed as, for
 $q=rt, p=\rho s, t, s\in \bS,$
\begin{eqnarray}\label{thanks2}
K(q,p)&=&\langle h_{q},h_{p}\rangle_{L^2(\bS)}\nonumber \\
&=&\int_{\bS} h_{q}(t')h_{p}(t')d\sigma(t')\nonumber \\
&=&\int_{\bS} \frac{1-r^2}{|q-t'|^d}\frac{1-\rho^2}{|p-t'|^d}d\sigma(t')\nonumber \\
&=& P_{\rho rt}(s)\\
&=& P_{r\rho s}(t)\nonumber .
\end{eqnarray}
The last equality is due to the relation $|r\rho t-s|^2=|r\rho s-t|^2.$   Now we prove the second last equality of the above equality chain. On one hand, the last integral is harmonic in $p,$ and when $\rho\to 1,$ the integral tends to $P_{q}(s)=P_{rs}(t).$ On the other hand, $P_{\rho rt}(s)=P_{r\rho s}(t)$ is also a harmonic function in $p,$ verified through the polar coordinate form of the Laplacian:
\[ \triangle_p=\frac{\partial^2}{\partial \rho^2}+\frac{d-1}{\rho}
\frac{\partial}{\partial \rho}+\frac{1}{\rho^2}\triangle_{\bS},\]
where $\triangle_{\bS}$ is the Beltrami-spherical Laplacian on the sphere. In fact,
\[ \triangle_p(P_{r\rho s}(t))=r^2\left(\frac{\partial^2}{\partial (r\rho)^2}+\frac{d-1}{(r\rho)}
\frac{\partial}{\partial (r\rho)}+\frac{1}{(r\rho)^2}\triangle_{\bS}\right)P_{(r\rho) s}(t)=0.\]
This harmonic function for $\rho\to 1$ has the same boundary limit function $P_{rs}(t).$ Due to the uniqueness of the solution of the harmonic boundary value problem the second last equality holds.

 The above is verification of the semi-group
 property of the Poisson kernel on the sphere. It provides computational conveniences.  With the $\mathcal{H}$-$H_K$ formulation the reproducing property
 of $K_q$ is automatic: For $u\in H_K=h^2({\bS}), q=rt,$
 \begin{eqnarray*}
 \langle u,K_q\rangle_{H_K}=\langle u(\cdot),P_{rt}(\cdot)\rangle_{L^2({\bS})}
 =u(rt)=u(q).\end{eqnarray*}
 To perform POAFD in $H_K=h^2(\bB)$ we need to prove the corresponding BVC.  For a general
   $K_{w}, w=\rho s, s\in \bS, 0<\rho < 1,$ its norm is computed
   \[  \|K_w\|^2_{H_K}=\langle K_w,K_w\rangle_{H_K}=K(w,w)=P_{\rho^2s}(s)
   =\frac{c_d(1+\rho^2)}{(1-\rho^2)^{d-1}}.\]
   The normalization of $K_q$ is, as denoted,
   \[ E_w=\frac{K_w}{\|K_w\|}=\frac{(1-\rho^2)^{(d-1)/2}}{\sqrt{c_d(1+\rho^2)}}K_w.\]

 We are to verify the BVC, that is, for $u$ being any function in $H_K,$
 \begin{eqnarray}\label{B1}
  \lim_{\bB\ni w\to \bS}|\langle u,E_w\rangle_{H_K}|=0.
  \end{eqnarray}
    We first recall that for $w=\rho s, s\in \bS,$
 \begin{eqnarray}\label{tk}
 \langle u,E_w\rangle_{H_K}=\frac{(1-\rho^2)^{(d-1)/2}}{\sqrt{c_d(1+\rho^2)}}u(w).
 \end{eqnarray}
Due to density of parameterized spherical Poisson kernels, verification of the BVC
for a general function $u\in H_K$
reduces to verifying the BVC for any but fixed parameterized
 spherical
 Poisson kernel $K_q, q=rt, t\in \bS.$ Thanks to the
 reproducing kernel expression (\ref{tk}) we have,
 \begin{eqnarray}\label{quantity1}E_w(q)=\langle K_q,E_w\rangle_{H_K}
 =\frac{(1-\rho^2)^{(d-1)/2}}{\sqrt{c_d(1+\rho^2)}}P_{r\rho t}(s).
 \end{eqnarray}
 When $\rho\to 1-,$ the quantities $P_{r\rho t}(s)$ are bounded uniformly in $t,s,$ and the bounds depend on the fixed $r<1.$ Since the factor in front tends to zero for $d\ge 2,$ the whole quantity tends to zero uniformly in $t,s.$  BVC is thus proved and POAFD performable.

 \begin{remark}
 The POAFD approximation obtained above amounts that for any positive integer
 $n$ and for function $u\in H_K^M$ there exists an
 $N$-combination of parameterized Poisson kernels
 that satisfies
 \[ \| u - \sum_{k=1}^N c_kP_{q_k}\|_{h^2({\bB})}\leq \frac{M}{\sqrt{N}}.\]
 Or, in terms of the boundary data in $L^2(\bS)$ it is
 \[ \| f(\cdot )- \sum_{k=1}^N c_kP_{q_k}(\cdot)\|_{L^2(\bS)}
 \leq \frac{M}{\sqrt{N}}.\]
 \end{remark}

\section{Experiments}
 Two experiments on sparse spherical Poisson kernel and heat kernel approximations are included. Below, the dot lines represent original functions, and the solid lines represent the approximation functions.

\begin{example} (Sparse Spherical Poisson Approximation) Let $f(q)$ be, as a toy example, the linear combination of three normalized spherical Poisson kernels on the $2$-sphere (d=3 with the formula (\ref{quantity1})) $E_{p_1}, E_{p_2}$ and $E_{p_3},$ reads $f(q)=\sum_{j=1}^3 c_jE_{p_j}(q),$ where $(c_1,c_2,c_3)=(0.8463,1.4105,0.0470), p_j=(\rho_j; \underline{s_j}), (\rho_1,\rho_2,\rho_3)=(0.4, 0.6, 0.8), \underline{s_j}=(\sin\phi_j \cos\theta_j,\sin\phi_j \sin\theta_j,\cos\phi_j),  \phi_j\in [0,\pi], \theta_j\in [0,2\pi),$\\ $(\phi_1,\phi_2,\phi_3)=\pi/5,\pi/2,=4\pi/5),
(\theta_1,\theta_2,\theta_3)=(\pi/5, 4\pi/5, 7\pi/5).$
Precisely, \begin{eqnarray}\label{extract}f(q)=\sum_{j=1}^{3}c_j \frac{1-\rho_j^2}{\sqrt{1+\rho_j^2}}\frac{1-(r\rho_j)^2}{|r\rho_j \underline{s}_j-\underline{t}|^3},\end{eqnarray}
where $q=(r;\underline{t}), \underline{t}=(\sin\phi \cos\theta,\sin\phi\sin\theta,\cos\phi), \phi\in [0,\pi], \theta\in [0,2\pi).$ Being only based on the boundary data extracted from (\ref{extract}), by doing POAFD with iterations 2,4,6,8 we obtain four POAFD expansions, with relative errors, respectively, 0.4310, 0.0237, 0.0022 and 0.3$\times 10^{-5}.$ The consecutive $8$ parameters $q_j=(h_j;\underline{w_j})$ with $ \underline{w_j}=(\sin\alpha_j \cos\beta_j,\sin\alpha \sin\beta_j,\cos\alpha_j),  \alpha_j\in [0,\pi], \beta_j\in [0,2\pi), j=1,\cdots,8,$ are\\
$(h_1,h_2,h_3,h_4,h_5,h_6,h_7,h_8)=
(0.4041,0.5714,0.6999,0.4518,0.5533,0.8207,0.4738,0.4042),$\\
$(\alpha_1,\alpha_2,\alpha_3,\alpha_4,\alpha_5,\alpha_6,\alpha_7,\alpha_8)
=(0.4518,0.7200,1.6200,0.6200,1.5200,1.6201,1.4201,1.6202),$ and $(\beta_1,\beta_2,\beta_3,\beta_4,\beta_5,\beta_6,\beta_7,\beta_8)=
(2.4200,0.3200,2.5200,0.5200,2.5201,4.4200,2.5202,0.6200).$
For showing the approximation efficiency of the recovering functions, we set $r=1, \theta=3.02,$ and the graphs are functions of $\phi$ varying in between $0$ and $\pi.$
\end{example}

\begin{figure}[H]
  \begin{minipage}[c]{0.22\textwidth}
  \centering
  \includegraphics[height=4.5cm,width=4cm]{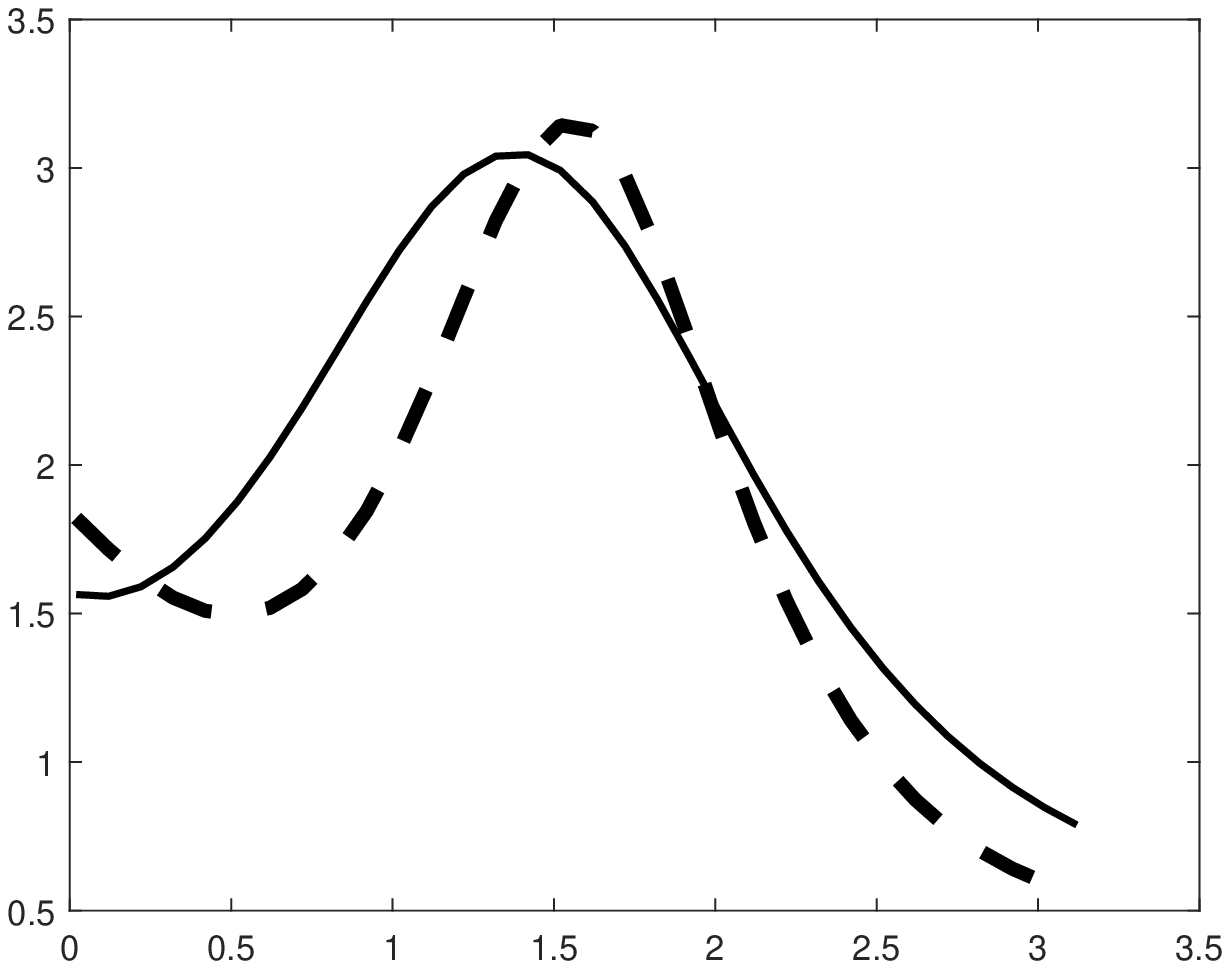}
  \caption*{{\tiny 2 POAFD iterations}}
  \end{minipage}
  \begin{minipage}[c]{0.22\textwidth}
  \centering
  \includegraphics[height=4.5cm,width=4cm]{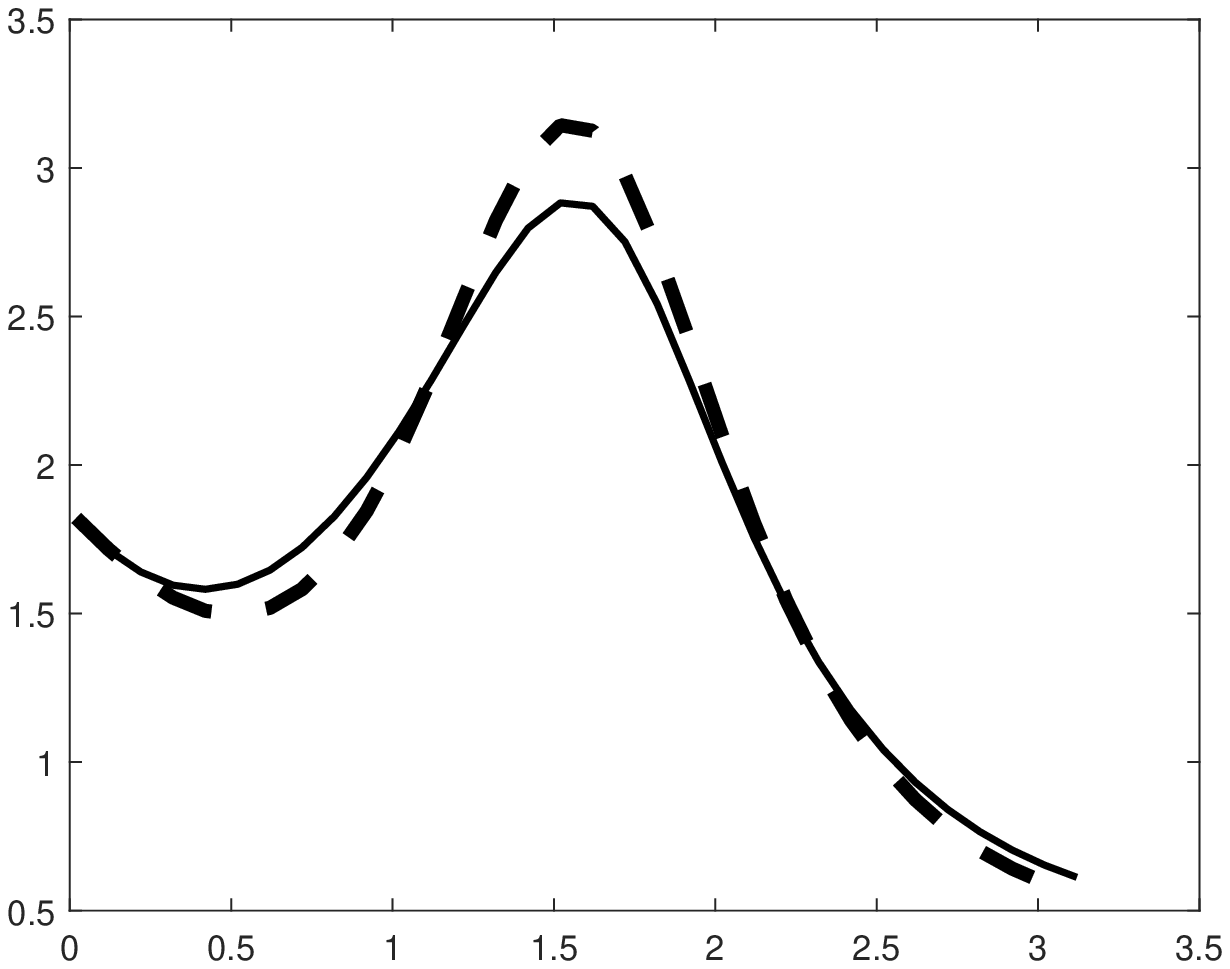}
  \caption*{{\tiny 4 POAFD iterations}}
  \end{minipage}
  \begin{minipage}[c]{0.22\textwidth}
  \centering
  \includegraphics[height=4.5cm,width=4cm]{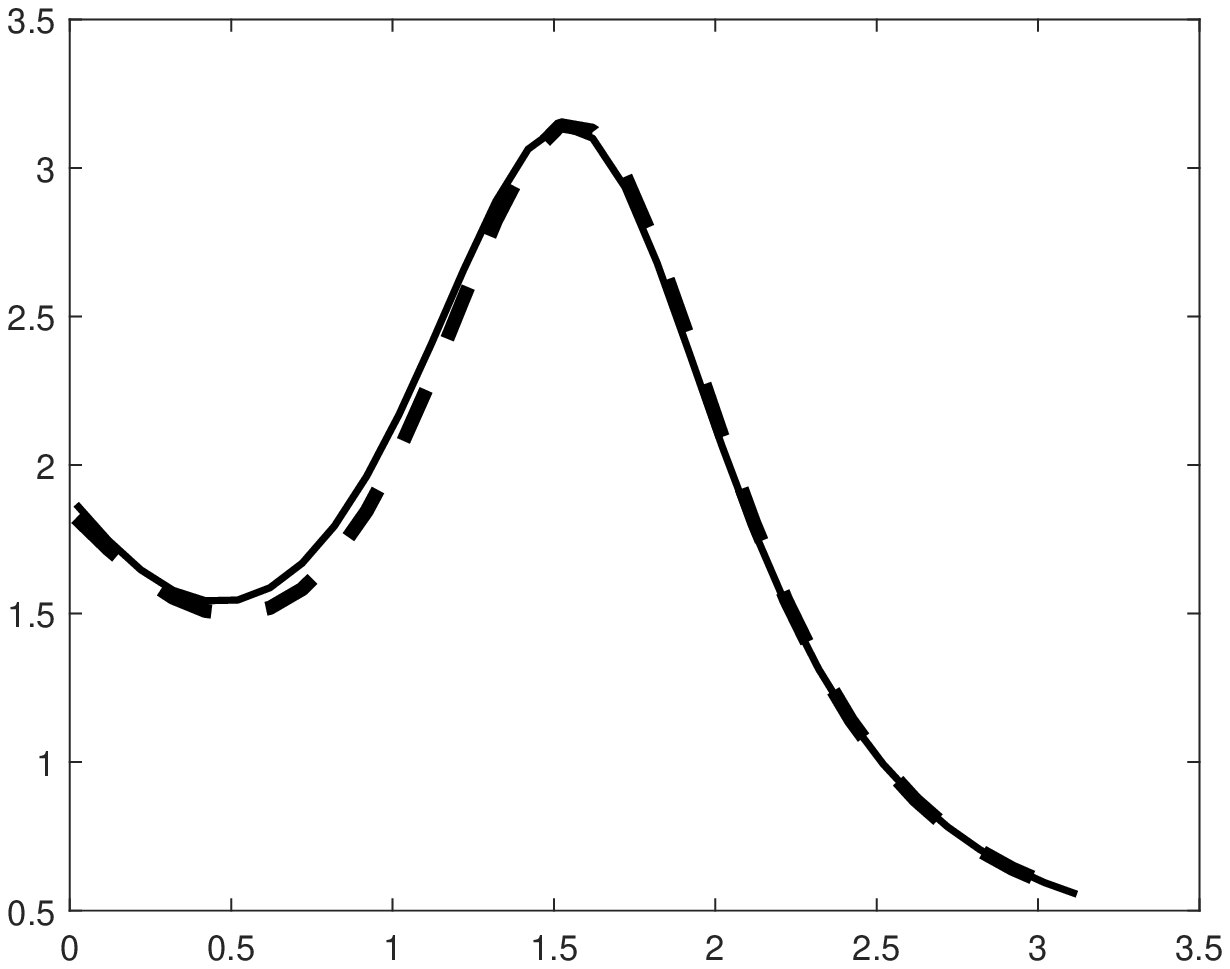}
  \caption*{{\tiny 6 POAFD iterations}}
  \end{minipage}
   \begin{minipage}[c]{0.22\textwidth}
  \centering
  \includegraphics[height=4.5cm,width=4cm]{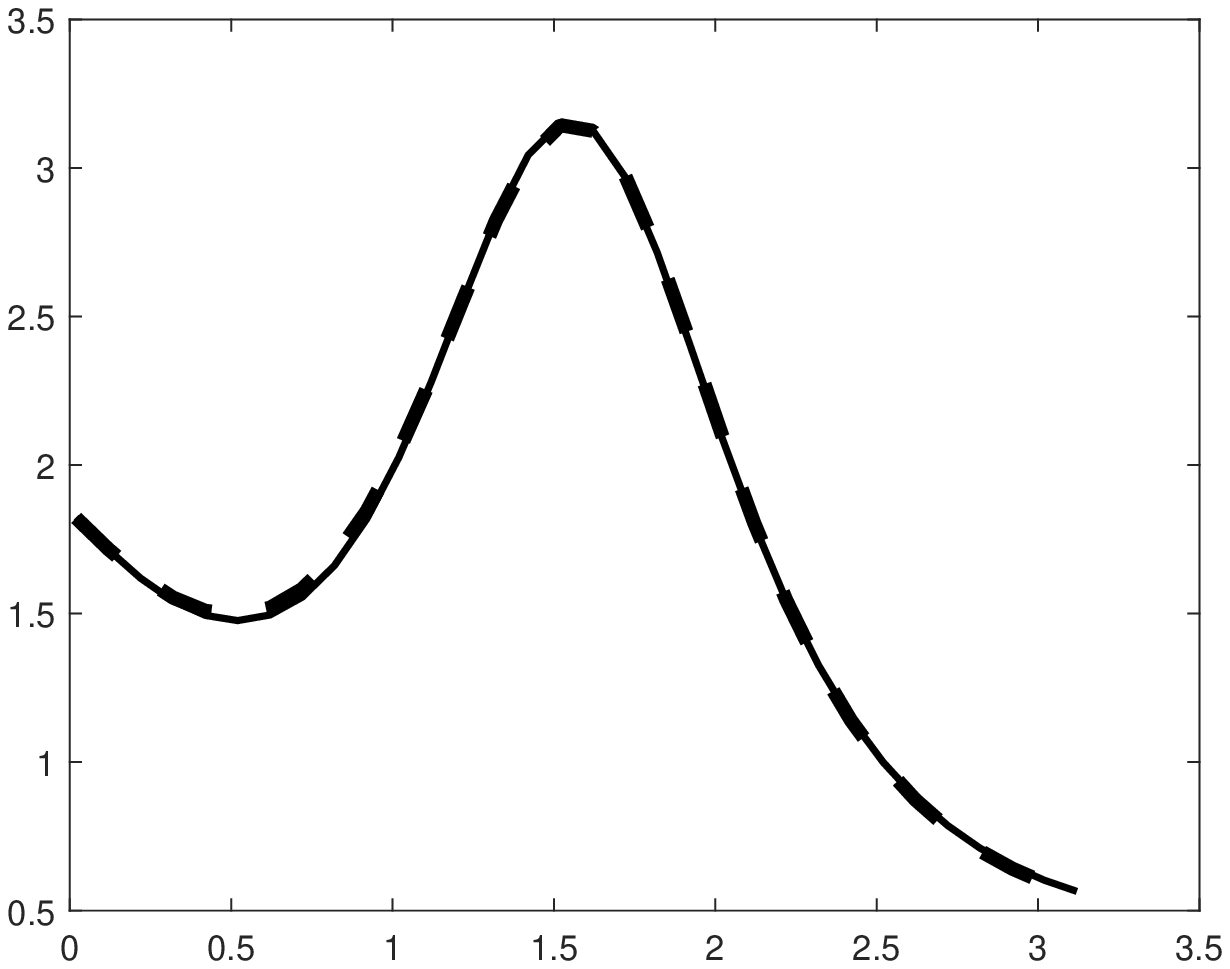}
  \caption*{{\tiny 8 POAFD iterations}}
  \end{minipage}
  \caption{Sparse Poisson kernel approximation on the unit Sphere}\label{Possion}
\end{figure}

\begin{example} (Sparse Heat Kernel Approximation) Let the signal to be expanded be given as
$$f(q)=\sum_{j=1}^{4}c_j\frac{\sqrt{8\pi t_j}}{4\pi(t+t_j)}e^{-\frac{|\underline{x}-\underline{y}_j|^2}{4(t+t_j)}},$$
where $q=(t,\underline{x}), t>0, \underline{x}=(x_1,x_2)\in {\R}^2,$ $p_j=(s_j,\underline{y_j}), j=1,2,3,4,\  (s_1,s_2,s_3,s_4)=(3,1,5,7),$ $\underline{y_1}=(-1,1), \underline{y_2}=(1,-5), \underline{y_3}=(2,6), \underline{y_4}=(-5,2),\ (c_1,c_2,c_3,c_4)=(0.05,0.5,0.01,1).$ The POAFD iteration numbers are 3,5,7 while the relative errors are 0.0190, 0.0087, 0.0002, respectively. By doing the POAFD the corresponding recovery parameters are $q_1=(10.0000;-4.1000,1.0000), q_2=(1.0000;1.200,-5.3000),
q_3=(7.6000;-5.1000,2.5000),
q_4=(4.4000;-1.6000,1.5000), q_5=(5.5000;2.5000,7.5000), q_6=(3.8000;-1.0000,1.4000),
q_7=(5.0000;2.0000,6.5000)$. The graphs of the recovering functions are for $t=0, x_2=-10.$
\end{example}

\begin{figure}[H]
  \begin{minipage}[c]{0.3\textwidth}
  \centering
  \includegraphics[height=5cm,width=4.5cm]{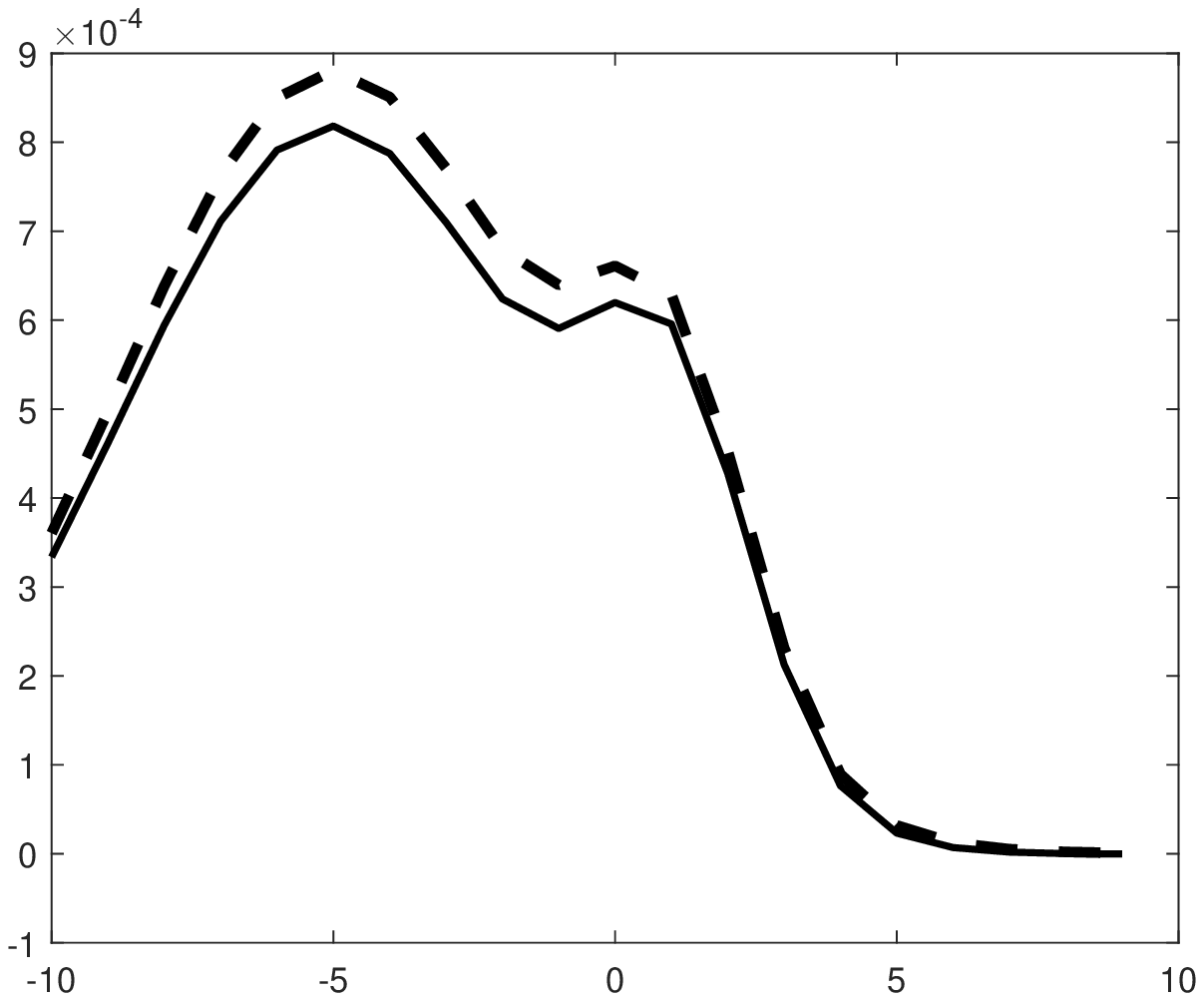}
  \caption*{{\tiny 3 POAFD iterations}}
  \end{minipage}
  \begin{minipage}[c]{0.3\textwidth}
  \centering
  \includegraphics[height=5cm,width=4.5cm]{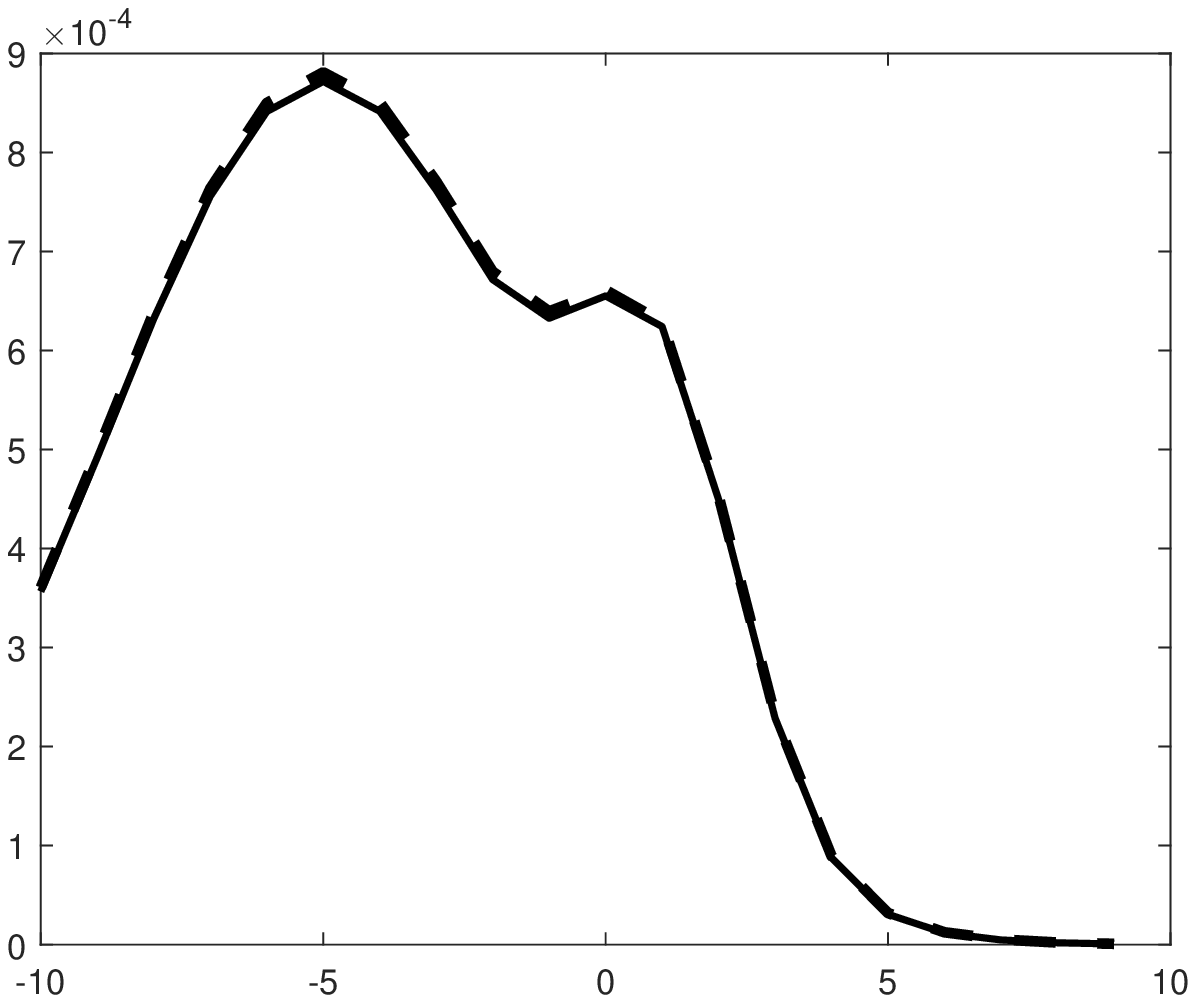}
  \caption*{{\tiny 5 POAFD iterations}}
  \end{minipage}
  \begin{minipage}[c]{0.3\textwidth}
  \centering
  \includegraphics[height=5cm,width=4.5cm]{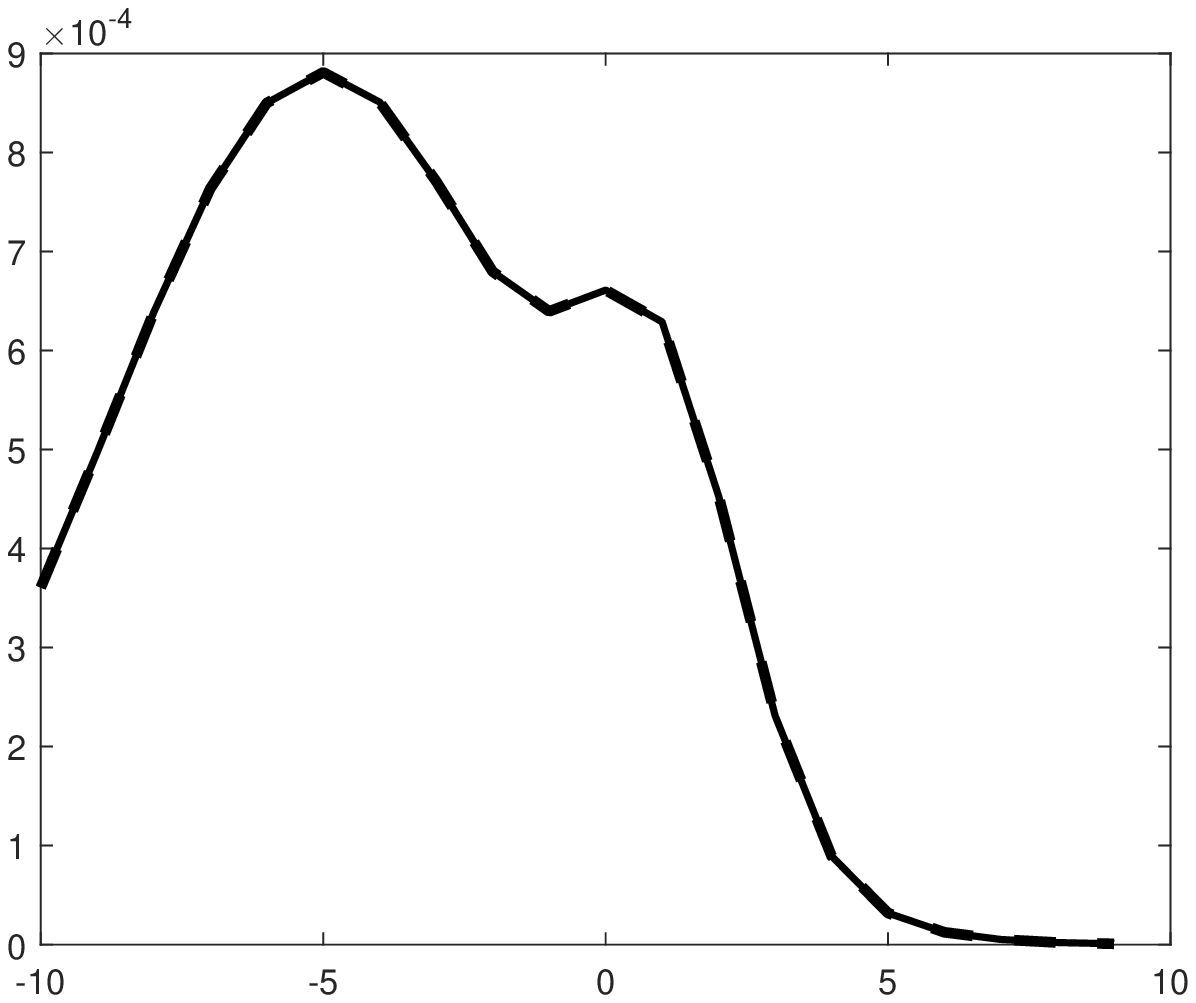}
  \caption*{{\tiny 7 POAFD iterations}}
  \end{minipage}
  \caption{Heat kernel approximation on the upper half space}
\end{figure}


\begin{thebibliography}{10}
\bibitem{ACQS1} D. Alpay, F. Colombo, T. Qian, I. Sabadini, {\it Adaptive orthonormal systems for matrix-valued functions}, Proceedings of the American Mathematical Society, 2017, 145(5)ㄩ2089每2106.

\bibitem{ACQS2} D. Alpay,  F. Colombo, T. Qian, and I. Sabadini, {\it Adaptative Decomposition: The Case of the Drury-Arveson Space}, Journal of Fourier Analysis and Applications, 2017, 23(6): 1426-1444.

\bibitem{Bara1} L. Baratchart, {\it Existence and generic properties of $L^2$ approximations for linear systems,} Math. Control Inform., {\bf 3}: 89-101.

\bibitem{Bara2} L. Baratchart, M. Cardelli, M. Olivi, {\it Identification and rational L2 approximation a gradient algorithm}, Automatica, 1991, 27(2): 413-417.

\bibitem{CQT} Q.-H. Chen, T. Qian, L.-H. Tan, {\it A Theory on Non-Constant Frequency Decompositions and Applications,} In: Advancements in Complex Analysis: From Theory to Practice, D. Breaz and M. Th. Rassias (Eds.), Springer, (to appear).

\bibitem{MZ}M. Gharavi-Alkhansari, T. S. Huang, {\it A fast orthogonal matching pursuit algorithm}, In Proceedings of the 1998 IEEE International Conference on Acoustics, Speech and Signal Processing, \textbf{3}:1389-1392.

\bibitem{LT} E. D. Livshitz, V. N. Temlyakov, {\it On convergence of weak greedy algorithms}, South Carolina Uiversity Columbia Department of Mathematics,  2000.

\bibitem{LZQ} Y. T. Li, L. M. Zhang, T. Qian, {\it 2D Partial Unwinding - A Novel Non-Linear Phase Decomposition of Images,} IEEE Transactions on Image Processing, 2019, DOI: 10.1109/TIP.2019.2914000.

\bibitem{MQ1} W.-X. Mai, T. Qian, {\em Rational Approximation in Hardy Spaces on Strips,} Complex Variables and Elliptic Equations, 2018, {\bf 63} : 1721-1738.

\bibitem{MQ2} W.-X. Mai, T. Qian, {\em Aveiro Method in Reproducing Kernel Hilbert Spaces under Complete Dictionary}, Mathematical Methods in the Applied Sciences, 2017, {\bf{40}}: 7240-7254.

\bibitem{Mi1} W. Mi, T. Qian, {\it Frequency-domain identification: An algorithm based on an adaptive rational orthogonal system,} {Automatica}, 2012, 48(6): 1154-1162.

\bibitem{Mi2} W. Mi, T. Qian, ``On backward shift algorithm for estimating poles of systems,'' {Automatica}, 2014, 50(6), 1603-1610.

\bibitem{MQW} W. Mi, T. Qian, F. Wan, {\it A Fast Adaptive Model Reduction Method Based on Takenaka-Malmquist Systems}, Systems and Control Letters, 2012, 61(1): 223每230.

\bibitem{MaZ}S. Mallat, Z. Zhang, {\it Matching pursuits with time-frequency dictionaries,} IEEE Trans. Signal Process, 1993, {\bf{41}}: 3397每3415.

    \bibitem{Q2020} T. Qian, {\it Reproducing Kernel Sparse Representations in Relation to Operator Equations,} Complex Anal. Oper. Theory 14 (2020), no. 2, 1每15.

\bibitem{Qian-cyclic}  T. Qian, {\it Cyclic AFD Algorithm for Best Rational}, Mathematical Methods in the Applied Sciences, 2014, 37(6): 846-859.
%
\bibitem{Q20} T. Qian, {\it A novel Fourier theory on non-linear phase and applications,} Advances in Mathematics (China), 2018, 47(3): 321-347 (in Chinese).

\bibitem{Qian1} T. Qian, {\it Reproducing Kernel Sparse Representations in Relation to Operator Equations}, Complex Analysis Operater Theory, 2020, 14(2):1每15.

\bibitem{Q2D} T. Qian, {\it Two-Dimensional Adaptive Fourier Decomposition}, Mathematical Methods in the Applied Sciences, 2016, 39(10) : 2431-2448.

\bibitem{QD1} W. Qu, P. Dang, {\it Rational approximation in a class of weighted Hardy spaces,} Complex Analysis and Operator Theory, 2019,13(4): 1827-1852.

\bibitem{QD2} W. Qu, P. Dang, {\it Reproducing kernel approximation in weighted Bergman spaces: Algorithm and applications,} Mathematical Methods in the Applied Sciences, 2019, 42(12): 4292-4304.

    \bibitem{QQD} W. Qu, T. Qian, D.G. Deng, {\it A sufficient condition for $n$-best approximation,} preprint.
%
\bibitem{QSW} T. Qian, W. Sproessig, J. X. Wang, {\it Adaptive Fourier decomposition of functions in quaternionic Hardy spaces,} Mathematical Methods in the Applied Sciences, 2012, 35(1): 43每64.

\bibitem{QWa} T.~Qian, Y.-B. Wang, {\it Adaptive Fourier series-a variation of greedy algorithm}, Advances in Computational Mathematics 2011,34~(3):279--293.

\bibitem{QWe} T. Qian, E. Wegert, {\it Optimal Approximation by Blaschke Forms,} Complex Variables and Elliptic Equations, 2013, 58(1): 123-133.

\bibitem{QWM} T. Qian, J. Z. Wang, W. X. Mai, {\it An Enhancement Algorithm for Cyclic Adaptive Fourier Decomposition}, Applied and Computational Harmonic Analysis, available online 19 January 2019.

\bibitem{Ruck1978} G. Ruckebusch, {\it Sur l'approximation rationnelle des filtres,} Report No 35 CMA Ecole Polytechnique, 1978.

\bibitem{SW} E. Stein, G. Weiss, {\it Introduction to Fourier Analysis in Euclidean Spaces,} Princeton University Press, 1970.

\bibitem{Saitoh1}S. Saitoh, Y. Sawano, {\it Theory of Reproducing Kernels and Applications, } Singapore: Springer, 2016.

\bibitem{Saitoh2} S. Saitoh, {\it Theory of Reproducing Kernels and its Applications,} Pitman Research Notes in Mathematics Series, vol. 189 (Longman Scientific and Technical, Harlow, 1988).

\bibitem{Te} V. Temlyakov, {\it Greedy Approximation,}  Cambridge Monographs an Applied and Computational Mathematics, 2011.

\bibitem{Walsh} J.~L. Walsh, {\it Interpolation and approximation by rational functions in the
  complex domain,} American Mathematical Soc. Publication, 1962.

\bibitem{WQLG} X. Y. Wang, T. Qian, I. T. Leong, Y. Gao, {\it Two-Dimensional Frequency-Domain System Identification}, IEEE Transactions on Automatic Control, 2019, DOI: 10.1109/TAC.2019.2913047.

\end{thebibliography}
\end{document}